\documentclass[10pt,a4paper,twoside]{article}
\usepackage{amsfonts}
\usepackage{amssymb}
\usepackage{amsthm}
\usepackage{latexsym}
\usepackage{amsmath}
\usepackage{amssymb}
\newtheorem{teo}{Theorem}[section]
\newtheorem{conj}[teo]{Conjecture}

\newtheorem{lema}[teo]{Lemma}
\newtheorem{cor}[teo]{Corollary}
\newtheorem{prop}[teo]{Proposition}
\newtheorem{defin}[teo]{Def\mbox{}inition}

\newtheorem{obs2}[teo]{Remark}

\newcounter{local}

\newenvironment{obs}{\begin{obs2}\rm}{\end{obs2}}

\newtheorem{no2}[teo]{Note}

\newenvironment{dem}{\begin{proof}[Proof]}{\end{proof}}

\newcommand{\Rr}{\mathcal{R}}
\newcommand{\C}{\mathbb{C}}
\newcommand{\A}{\mathbb{A}}
\newcommand{\Q}{\mathbb{Q}}
\newcommand{\Z}{\mathbb{Z}}
\newcommand{\Y}{\mathcal{O}}

\newcommand{\f}{\mathfrak{f}}
\newcommand{\p}{\mathfrak{p}}

\newcommand{\Ff}{\mathcal{F}}

\newcommand{\ilim}[1]{\displaystyle{\lim_{\stackrel{\longrightarrow}{#1}}}}
\begin{document}
\begin{center}
\LARGE{On Jannsen's conjecture for Hecke characters of imaginary quadratic fields}\\
\vspace{0.5cm}
\Large{Francesc Bars \footnote{Work partially supported by BFM2003-06092}}\\

\end{center}

\begin{center}
\begin{small}
 \begin{abstract}
 We present a collection of results on a conjecture of Jannsen about
the $p$-adic realizations associated to Hecke characters over an
imaginary quadratic field $K$ of class number 1.

The conjecture is easy to check for Galois groups purely of local
type (\S 1). In \S2 we define the $p$-adic realizations associated
to Hecke characters over $K$. We prove the conjecture under a
geometric regularity condition for the imaginary quadratic field
$K$ at $p$, which is related to the property that a global Galois
group is purely of local type. Without this regularity assumption
at $p$, we present a review of the known situations in the
critical case \S3 and in the non-critical case \S4 for these
realizations. We relate the conjecture to the non-vanishing of
some concrete non-critical values of the associated $p$-adic
$L$-function of the Hecke character.

Finally, in \S 5 we prove that the conjecture follows from a
general conjecture on Iwasawa theory for almost all Tate twists.
\end{abstract}
\end{small}
\end{center}

\begin{section}{The Jannsen conjecture on local type Galois groups}
Jannsen's conjecture \cite{Ja} predicts the vanishing of a second
Galois cohomology group for the $p$-adic realization of almost all
Tate twists of a pure Chow motive. It also specifies the Tate
twists where this cohomology group could not vanish. Without this
specification the conjecture is a generalization of the classical
weak Leopoldt conjecture.

 We refer to Jannsen's original
paper \cite{Ja} and Perrin-Riou's paper \cite[Appendix B]{PR1} for
the relations with other conjectures and for general results.

Let $F$ be a number field with algebraic closure $\overline{F}$.
Let $X$ be a smooth, projective variety of pure dimension $d$ over
$F$. Let $p$ be a prime number, and $S$ a finite set of places of
$F$, containing all places above $\infty$ and $p$, and all primes
where $X$ has bad reduction. Let $G_S$ be the Galois group over
$F$ of the maximal $S$-ramified (unramified outside $S$) extension
of $F$, that we call $F_S$.
\begin{conj}[Jannsen] If $\overline{X}=X\times_{F}\overline{F}$, then
$$H^2(G_S,H^i_{et}(\overline{X},\Q_p(n)))=0\ if \ \left\{\begin{array}{l@{\quad \quad}l}
a)\ i+1<n,& or\\
b)\ i+1>2n.\\
\end{array}\right.$$
\end{conj}
This conjecture can be verified also from the \'etale cohomology
with $\Z_p$ or $\Q_p/\Z_p$ coefficients.
\begin{lema}[lemma 1 \cite{Ja}] \label{lema12} The following statements are equivalent:
\begin{enumerate}
\item $H^2(G_S,H^i_{et}(\overline{X},\Q_p(n)))=0$. \item
$H^2(G_S,H^i_{et}(\overline{X},\Z_p(n)))$ is finite. \item
$H^2(G_S,H^i_{et}(\overline{X},\Q_p/\Z_p(n)))$ is finite. \item
(if $p\neq 2$ or $F$ is totally imaginary)
$H^2(G_S,\tilde{H}^i_{et}(\overline{X},\Q_p/\Z_p(n)))=0$ where
$\tilde{H}^i_{et}(\overline{X},\Q_p/\Z_p(n))$ is the $p$-divisible
part of the group $H^i_{et}(\overline{X},\Q_p/\Z_p(n))$.
\end{enumerate}

\end{lema}
There is an analogous conjecture if one replaces $G_S$ by
$\operatorname{Gal}(\overline{F}_{\p}/F_{\p})$, the absolute
Galois group of the local field $F_{\p}$, where $F_{\p}$ is the
completion at $\p$ of $F$.

\begin{conj}[Jannsen] $H^2(\operatorname{Gal}(\overline{F_{\p}}/F_{\p}),H^i_{et}(\overline{X},\Q_l(n)))=0$ if
$i+1<n$  or $i+1>2n$.
\end{conj}

Lemma \ref{lema12} is also valid in this situation, i.e. replacing
$G_S$ by $\operatorname{Gal}(\overline{F_{\p}}/F_{\p})$.

Denote by $G_S(p)$ the Galois group
$\operatorname{Gal}(F_S(p)/F)$, where $F_S(p)$ means the maximal
$p$-extension of $F$ inside $F_S$. Let $S_p$ be the set of primes
of $F$ above $p$.

Let us denote by $\mathcal{G}_v$ the Galois group of the maximal
$p$-extension of the local field $F_v$ (with $v\in S$) over $F_v$;
we write also $\mathcal{G}_v=\operatorname{Gal}(F_v(p)/F_v)$. We
have the natural surjective map
$$i_{v,S}:\mathcal{G}_v\rightarrow G_v$$
where $G_v$ denotes the $p$-part of the decomposition group at
$F_v$ of $G_S$. Let $\varphi_v:G_v\rightarrow G_S(p)$ be the
natural inclusion.
\begin{defin} The Galois group $G_S(p)$ is purely of local type iff the map $\varphi_v\circ i_{v,S}$
is an isomorphism.
\end{defin}
It is easy to see that this definition is equivalent to the one
given in \cite{W2}(\cite[B.2.6]{Ba}).
\begin{lema}\label{Neu} Let $M(j)$ be a $p$-primary divisible $\mathcal{G}_v=\operatorname{Gal}({F_{v}}(p)/F_{v})$-module
of cofinite type. If $\mu_p\notin F_{v}$ then
$H^2(\mathcal{G}_v,M(j))=0$, where $\mu_p$ denotes a primitive
$p$-root of unity. If $\mu_p\in F_{v}$, then
$H^2(\mathcal{G}_v,M(j))=M(j-1)_{\mathcal{G}_v}$.

Finally, let $M(j)$ be any $p$-primary divisible
$\operatorname{Gal}(\overline{F}_{v}/F_{v})$-module of cofinite
type. Then
$H^2(\operatorname{Gal}(\overline{F}_{v}/F_{v}),M(j))=M(j-1)_{\operatorname{Gal}(\overline{F}_{v}/F_{v})}$.
\end{lema}
\begin{dem} If $\mu_p\notin F_{v}$ then $\mathcal{G}_v$ is free, hence
$H^2(\mathcal{G}_v,M(j))=0$. Otherwise it is a Poincar\'e group of
dimension two with dualizing module $\Q_p/\Z_p(1)$. Using local
Tate duality we get
$H^2(\mathcal{G}_v,M(j))=M(j-1)_{\mathcal{G}_v}$.

The last statement follows from local Tate duality.
\end{dem}
\begin{lema}\label{lema1.6}Let $S'$ be a subset of $S$ which contains $S_p$, the places of $F$ above $p$.
Suppose that the Galois group $\operatorname{Gal}(F_{S'}(p)/F)$ is
purely of local type for a place $w\in S'$ of $F$. Furthermore let
$M$ be a $p$-primary divisible
$\operatorname{Gal}(F_{S'}(p)/F)$-module of cofinite type such
that $M(j-1)_{\operatorname{Gal}(\overline{F}_v/F_v)}=0$ for all
$v\in S\setminus(S'\setminus\{w\})$. Then
$$H^2(G_S,M(j)=0.$$
\end{lema}
\begin{dem} Let us consider part of Soul\'e's exact sequence recalled in
\cite[\S3,(13)]{Ja}:

$$H^2(G_{S'},M(j))\rightarrow H^2(G_S,M(j))\rightarrow
\oplus_{v\in S\setminus
S'}M(j-1)_{\operatorname{Gal}(\overline{F}_v/F_v)}\rightarrow 0,$$
where the last term of this sequence is zero by hypothesis.
Observe that
$$H^2(G_{S'},M(j))=H^2(G_{S'}(p),M(j))$$ by a result of Neumann \cite[Theorem 1]{Neum}.
Since $G_{S'}(p)$ is purely of local type,
$$H^2(G_{S'}(p),M(j))=H^2(\mathcal{G}_w,M(j)).$$
If $\mu_p\notin F_{w}$ then $H^2(\mathcal{G}_w,M(j))=0$. Otherwise
$H^2(\mathcal{G}_w,M(j))=M(j-1)_{\operatorname{Gal}({F}_{w}(p)/F_{w})}$.
But, in our case, the action of
$\operatorname{Gal}(\overline{F}_w/F_w)$ factors through the
pro-$p$-quotient $\mathcal{G}_w$, and so
$M(j-1)_{\operatorname{Gal}({F}_{w}(p)/F_{w})}=0$ by hypothesis.
\end{dem}

\begin{obs2}\label{janloctrue} If $X$ has good reduction at $\p$ and $\p\nmid p$, the vanishing of the Galois group $H^2(F_{\p},H^i_{et}(\overline{X},\Q_p(j)))$
with $i\neq 2(j-1)$ follows from lemma \ref{Neu} and weights
arguments (see \cite[lemma 3]{Ja}). As usual $H^i(F_v,-)$ means
$H^i(\operatorname{Gal}(\overline{F_v}/F_v),-)$. This result can
also be proved when $X$ has potentially good reduction at $\p$
\cite[lemma 12]{Ja}.

When $\p | p$, Soul\'e proves the vanishing of this group for
$j>i+1$ or $j<1$ \cite[corollary 5,lemma 11]{Ja}. Moreover if $\p$
is unramified in $F/\Q$ and $X$ has good reduction at $\p$ then
$H^2(F_{\p},H^i_{et}(\overline{X},\Q_p(j)))=0$  if $i\neq 2(j-1)$
\cite[corollary 6]{Ja}. In this situation if $i-2j\neq2$ and
$i-2j\neq 0$ we obtain:
$$\operatorname{dim}_{\Q_p}H^1(F_{\p},H^i_{et}(\overline{X},\Q_p(j)))=-\chi(H^i_{et}(\overline{X},\Q_p(j)))=$$
$$[F_{\p}:\Q_p]\operatorname{dim}_{\Q_p}(H^i_{et}(\overline{X},\Q_p(j))),$$
where $\chi$ is the Euler characteristic in local Galois
cohomology. The last equality follows from \cite[7.3.8]{NWS}.
\end{obs2}

\end{section}
\begin{section}{The geometric regularity condition.}

We keep once and for all the notations of \S 1. Let
 $\p$ be a prime of $F$ such that $\p\mid p$ and the
inertia group at $\p$ acts trivially on $M$, or on $\overline{X}$.

 Let $E$ be a fixed elliptic curve
 defined over an imaginary quadratic field $K$,
and with CM by $\Y_K$, the ring of integers of $K$. This
hypothesis implies that $cl(K)=1$. Associated to this elliptic
curve there is a Hecke character $\varphi$ of the imaginary
quadratic field $K$ with conductor $\f$, an ideal of $\Y_K$ which
coincides with the conductor of the elliptic curve $E$.

  Consider the category of Chow motives $\mathcal{M}_{\Q}(K)$ over
$K$ with morphisms induced by graded correspondences in Chow
theory tensored with $\Q$. Then, the motive of the elliptic curve
$E$ has a canonical decomposition $h(E)_{\Q}=h^0(E)_{\Q}\oplus
h^1(E)_{\Q}\oplus h^2(E)_{\Q}$. The motive $h^1(E)_{\Q}$ has a
multiplication by $K$ \cite{De}. For $w$ a strictly positive
integer, let us consider the motive $\otimes^w_{\Q} h^1(E)_{\Q}$,
which has multiplication by $T_w:=\otimes^w_{\Q}K$. Observe that
$T_w$ has a decomposition $\prod_{\theta}T_{\theta}$ as a product
of fields $T_{\theta}$, where $\theta$ runs through the
$Aut(\C)$-orbits of $\Upsilon^w=\operatorname{Hom}(T_w,\C)$, where
$\Upsilon=\operatorname{Hom}(K,\C)$. This decomposition defines
some idempotents $e_{\theta}$ and gives also a decomposition of
the motive and its realizations. Let us fix once and for all an
immersion $\lambda:K\rightarrow\C$ as in \cite[p.135]{De}.

The $L$-function associated to the motive $e_{\theta}(\otimes^w
h^1(E)_{\Q})$ corresponds to the $L$-function associated to the
CM-character
$\psi_{\theta}=e_{\theta}(\otimes^w\varphi):\A_K^*\rightarrow K^*$
\cite[\S 1.3.1]{De} (for equivalent definitions of Hecke
characters over an imaginary quadratic field, see \cite[\S
2.2]{G}). With respect to the fixed embedding $\lambda$, this
CM-character corresponds to $\varphi^a\overline{\varphi}^b$, where
$a,b\geq0$ are integers such that $w=a+b$. The pair $(a,b)$ is the
infinite type of $\psi_{\theta}$. We note that there are different
$\theta$ with the same infinite type. Every $\theta$ gives two
elements of $\Upsilon^w$, one given by the infinite type
$\vartheta\in\theta\cap \operatorname{Hom}_{K}(T_w,\C)$ and the
other one coming from the other embedding. When $\theta$
corresponds to $\vartheta=(\lambda_1,\ldots,\lambda_w)\in
\operatorname{Hom}_{K}(T_w,\C)$ with type $(a,b)$, we denote by
$\overline{\theta}$ the one that corresponds to the element
$\overline{\vartheta}=(\overline{\lambda_1},\ldots,\overline{\lambda_w})\in
\operatorname{Hom}_{K}(T_w,\C)$, where $\lambda_i$ are
$\Q$-immersions of $K$ into $\C$ and the bar denotes complex
conjugation. Observe that $\overline{\theta}$ has type $(b,a)$.

Let us denote by $\otimes^wh^1(E)$ the integer Chow motive
(similar to $\otimes^wh^1(E)_{\Q}$ but without tensoring by $\Q$
the correspondences). The Chow motive
$$M_{\theta}:=e_{\theta}(\otimes^wh^1(E)),$$
is a Chow motive with coefficients in $\mathcal{O}_K[1/D_K]$ where
$D_K$ is the discriminant of $K$, $e_{\theta}\in
\mathcal{O}_K[1/D_K]$ when $w>1$.

The $p$-adic realization of the motive $M_{\theta}(n)$:
$$M_{\theta\Q_p}(n):=H_{et}^w(M_{\theta}\times_{K}\overline{K},\Q_p(n))\cong
e_{\theta}(\otimes^w (T_pE\otimes\Q_p))(n-w),$$ is unramified
outside any set $S$ which contains the finite primes of $K$ which
divide $p\f_{\theta}$, where $\f_{\theta}$ is the conductor of
$\psi_{\theta}$. We have also a natural lattice associated to it,
corresponding for $p\nmid D_K$ to
$H_{\acute{e}t}^w(M_{\theta}\times_{K}\overline{K},\Z_p(n))\cong
e_{\theta}(\otimes^w T_pE)(n-w)$ and we denote it by
$$M_{\theta\Z_p}(n):=H_{\acute{e}t}^w(M_{\theta}\times_{K}\overline{K},\Z_p(n)).$$
We impose once and for all that $p\nmid D_K$.

Since $e_{\theta}(\otimes^wh^1(E))\subseteq h^w(E^w)$ Jannsen's
conjecture for $h^w(E^w)$ implies the following conjecture for
Hecke characters.
\begin{conj} \label{conj2.1} Let $p$ be a prime such that $E$ has good reduction
at the primes over $p$ in $K$. Let $S$ be a set of primes of $K$
which contains the primes of $K$ over $p$ and the primes of $\f$.
Then, $$H^2(\operatorname{Gal}(K_S/K),M_{\theta\Q_p}(n))=0,$$ for
$n>w+1$ or $w+1>2n$.
\end{conj}

 Consider $\underline{F}\subseteq K(E[p])$ with $K\subseteq \underline{F}$. We suppose also $p>3$.
 Denote by $\Ff=K(E[p])$. We know that
$E\times_K\Ff$ has good reduction everywhere.

We introduce the notion of regularity for $\p$ at $\underline{F}$,
which turns out to be closely related to the condition that a
global Galois group is purely of local type.

Let us assume in this section once for all that $p$ splits in $K$
with $p=\p\p^*$ where $\p,\p^*$ are different primes of $K$.

\begin{defin}\label{defin2.2} Suppose that $p>3$, $(p)=\p\p^*$ in $K$, $\p\neq\p^*$ and $E$ has good
reduction at $\p$ and $\p^*$. Let $S_{\p}$ be the set of primes of
$\underline{F}$ that divide $\p$. The prime $\p$ is called regular
for $E$ and $\underline{F}$ if $\underline{F}_{S_{\p}}(p)$ is a
$\Z_p$-extension of $\underline{F}$. We say that $p$ is regular
for $E$ and $\underline{F}$ if $\p$ and $\p^*$ are regular for $E$
and $\underline{F}$.
\end{defin}
\begin{obs2} The above regularity
condition admits an analogue of Kummer's criterion as in the
classical notion of regularity at $p$
\cite[B.1.2,1.7,2.16-17]{Ba}. One of the equivalent notions of
regularity is related to a critical value of the $L$-function
associated to some concrete Hecke characters of the form
$\varphi^a\overline{\varphi}^b$ \cite{Win} (\cite[B.1.7]{Ba}).
\end{obs2}

\begin{prop}[Wingberg \cite{Win}]\label{gua:teo} Let $\underline{F}$ be a number field between
$K$ and $\mathcal{F}$. The prime $\p$ is regular for $E$ and
$\underline{F}$ if and only if
$\operatorname{Gal}(\underline{F}_{S_p}(p)/\underline{F})$ is
purely of local type with respect to $\p^*$, where $S_p$ denotes
the places of $\underline{F}$ above $p$.
\end{prop}
\begin{obs2} Thus, for this notion of regular primes the Grunwald-Wang
theorem holds: the maximal $p$-extension of $\underline{F}_{\p^*}$
coincides with the completion at $\p^*$ of the maximal
$p$-extension of $\underline{F}_{S_p}$ (use \cite[corollary
B.2.6]{Ba} with proposition \ref{gua:teo}).
\end{obs2}

Let us write a consequence of proposition \ref{gua:teo} (the case
$w=1$ is a result of Wingberg \cite{Win}).
\begin{cor}\label{prewin} Let $p$ be a regular prime for $E$ and $\Ff=K(E[p])$.
 Let $\underline{F}$ be an extension of $K$ inside $\Ff$. Denote by $S^*$ a
finite set of primes of $\underline{F}$ containing the primes
above $p$ and those where $E\times_K \underline{F}$ has bad
reduction. Then,
$$H^2(\operatorname{Gal}(\underline{F}_{S^*}/\underline{F}),H^w(M_{\theta},\Q_p/\Z_p(j+w)))=0$$
for any integer $w=a+b$ with $w+2(j-1)\neq 0$.
\end{cor}
\begin{dem} By Kummer theory, we have that
$$H^1_{et}(\overline{E},\Q_p/\Z_p(1))=E[p^{\infty}]=E[\p^{\infty}]\oplus E[(\p^*)^{\infty}]$$ where $E[a^{\infty}]=\ilim{n}E[a^n]$ is the inductive
limit of the $a^n$-torsion points of the elliptic curve $E$. Thus
$M(w+j):=H^w(M_{\theta}\otimes_K\overline{K},\Q_p/\Z_p(w+j))\cong
e_{\theta}(\otimes^wE[p^{\infty}])(j)$. Since $K\subseteq
\underline{F}\subseteq\Ff$, we need only to prove
 $$H^2(\operatorname{Gal}(\Ff_{S^*}/\underline{F}),M(w+j))=0,$$
because $\underline{F}\subseteq\Ff$ is unramified outside
 $S^*$ \cite[II.1.8]{dS}. Since $\Ff/\underline{F}$ is prime
 to $p$ it is enough to show
 $$H^2(\operatorname{Gal}(\Ff_{S^*}/\Ff),M(w+j))=0.$$

Now, $M(w+j)$ is a $\operatorname{Gal}(\Ff_{S_p}(p)/\Ff)$-module
because the Galois action on $M(w+j)$ factors through
$\operatorname{Gal}(\Ff(E[p^{\infty}])/\Ff)$ and $E\times_K\Ff$
has good reduction \cite[1.9,(i),(ii)]{dS}.

 Using lemma \ref{lema1.6} it is enough to prove the vanishing of some coinvariant modules in
 local Galois groups.

  Since $E\times_K\Ff$ has good reduction
everywhere, in particular in the places of $S^*\setminus S_p$, the
same is true for $E^w=E\times_K E\times_K\ldots\times_K E$ over
the field $\Ff$. Then, by the proved Weil conjectures on
$H^w(\overline{E}^w,\Q_p(n))$, the Frobenius at $v$ over $\Ff$
does not act as identity in any subspace of the above cohomology
group for $w-2n\neq 0$. By local Tate duality
$M(w+j-1)_{\operatorname{Gal}(\overline{\Ff_v}/\Ff_v)}$ is dual to
$M_{\theta\Z_p}^{\operatorname{Gal}(\overline{\Ff_v}/\Ff_v)}$ and
this module vanishes for $w-2(w+j-1)\neq0$.

By hypothesis $\p^*$ is regular in $\Ff$, hence (theorem
\ref{gua:teo})
$$\operatorname{Gal}(\Ff_{S_p}(p)/\Ff)\cong \operatorname{Gal}(\Ff_{\p}(p)/\Ff_{\p}).$$  Thus, it
remains only to show that
$$M(w+j-1)_{\operatorname{Gal}(\Ff_{\p}(p)/\Ff_{\p})}=M(w+j-1)_{\operatorname{Gal}(\overline{\Ff_{\p}}/\Ff_{\p})}=0.$$
The $\operatorname{Gal}(\overline{\Ff_{\p}}/\Ff_{\p})$-action on
$M(w+j-1)$ factors though the pro-$p$-quotient by the regularity
and the condition $\mu_p\subseteq\Ff_{\p}$. By \cite[proof of
Lemma 3.2.5]{G} these coinvariants are always zero unless
$a-(j-1+w)=0$ and $b-(j-1+w)=0$. If both are zero, we have
$0=a-(j-1+w)+b-(j-1+w)=w-2j+2-2w=-(w+2j-2)$ which is impossible by
our hypothesis.
\end{dem}
\begin{obs2} The above result, and corollaries \ref{meujan2} and \ref{dosdeu} below, are true assuming
only that $S^*$ is a finite set of places which contains the
primes above $p$ and the primes of $\underline{F}$ such that the
inertia acts non-trivially in $M_{\theta}\times_K\underline{F}$
(use standard arguments like in the proof of \cite[2.2.16]{K}).
\end{obs2}
\begin{cor}[Jannsen's conjecture under regularity]\label{meujan2} Let $p$ be a regular prime for $E$ and $\Ff=K(E[p])$. Let $\underline{F}$ be
an extension of $K$ inside $\Ff$. Let $S^*$ be a finite set of
primes containing the primes of $\underline{F}$ that are over $p$
and the primes where the elliptic curve $E\times_K\underline{F}$
has bad reduction. Let $M_{\theta}$ be $e_{\theta}(\otimes^w
h^1(E))$. Then
$$H^2(\operatorname{Gal}(\underline{F}_{S^*}/\underline{F}),M_{\theta\Q_p}(n))=0$$
if $w+2(n-w-1)\neq 0$.
\end{cor}
\begin{dem}
Is a direct consequence of the above corollary \ref{prewin} and
lemma \ref{lema12}.
\end{dem}

\begin{cor}\label{meujan} Let $p$ be a regular prime for $E$ and $\Ff=K(E[p])$. Let $\underline{F}$ be
an extension of $K$ inside $\Ff$. Let $S^*$ be a finite set of
primes of $\underline{F}$ containing the primes above $p$ and the
primes where the elliptic curve $E\times_K\underline{F}$ has bad
reduction. Fix an integer $w\ge1$. Then
$$H^2(\operatorname{Gal}(\underline{F}_{S^*}/\underline{F}),(\otimes^wh^1(E)(1))_{et}(j))=0$$
for any integer $j$ such that $w+2(j-1)\neq0$, where
$(\otimes^wh^1(E)(1))_{et}(j))$ denotes the \'etale cohomology
group $(\otimes^w H^1_{et}(E\times_K\overline{K},\Q_p(1)))(j)$
which is the $p$-adic realization of the motive
$(\otimes^wh^1(E))(w+j)$.
\end{cor}
\begin{dem}
The integer idempotents $e_{\theta}$ give us a decomposition of
the above motive
$$(\otimes^wh^1(E)(1))(j)=\prod_{\theta}e_{\theta}((\otimes^wh^1(E)(1))(j)),$$
thus a decomposition of the $p$-adic realization
$$H^2(\operatorname{Gal}(\underline{F}_{S^*}/\underline{F}),(\otimes^wh^1(E)(1))_{et}(j))=$$
$$\oplus_{\theta}H^2(\operatorname{Gal}(\underline{F}_{S^*}/\underline{F}),e_{\theta}(\otimes^wh^1(E)(1))_{et}(j))$$
which is zero by Corollary \ref{meujan2}.
\end{dem}
\begin{cor}\label{dosdeu} Let $p$ be a regular prime for $E$ and $\Ff$ and let be $w-2n\neq0$ and
$2n-w-2\neq 0$. Let $S$ be a finite set of primes of $K$
containing the primes over $p$ and those where $E$ has bad
reduction. Then
$$\operatorname{dim}_{\Q_p}H^1(\operatorname{Gal}(K_S/K),M_{\theta\Q_p}(n))+\operatorname{dim}_{\Q_p}H^1(\operatorname{Gal}(K_S/K),M_{\overline{\theta}\Q_p}(w+1-n))=$$
$$=[K:\Q]\operatorname{dim}_{\Q_p}M_{\theta\Q_p}=4.$$
Moreover,
$$\operatorname{dim}_{\Q_p}\!H^1(\operatorname{Gal}(K_S/K),M_{\theta\Q_p}(w+l+1))\!=2=\operatorname{dim}_{\Q_p}\!H^1(\operatorname{Gal}(K_S/K),
M_{\overline{\theta}\Q_p}(-l))$$ for $-l\leq
\operatorname{min}(a,b)$ if $a\neq b$(and $-l<a$ if $a=b$).
\end{cor}
\begin{dem} By the conditions on the weights
and the regularity assumption we obtain
$H^0(G_S,M_{\theta\Q_p}(m))=0$ and $H^2(G_S,M_{\theta\Q_p}(m))=0$,
with $m=n$ and $m=w+1-n$ (and the same statements for
$\overline{\theta}$ instead of $\theta$). Thus, the first equality
holds by \cite[corollary 1, proof of lemma 2]{Ja}. Observe that
$M_{\theta\Q_p}\cong e_{\theta}(\otimes^wT_pE(-1))\otimes\Q_p$ has
$\Q_p$-rank equal to 2 (see for example \cite[\S 2]{Ba2});
therefore the second equality follows. The last statement follows
from the Beilinson conjecture for these characters (proved by
Deninger \cite{De}) and \cite[Lemma 2]{Ja} at the twist $w+l+1$.
\end{dem}
Observe $M_{\theta}(n)$ is a submotive of $h^w(E^w)(n)$, therefore
using Soul\'e's result (see remark \ref{janloctrue}) we obtain,
without the regularity assumption:
\begin{cor}\label{dosonze} Let $\p$ be a prime of $K$ such that
$\p|p$. Suppose $w-2n\neq-2$ if $\p$ is unramified, otherwise
$n>w+1$ or $n<1$. Then,
$$H^2(K_{\p},M_{\theta\Q_p}(n)))=0.$$ Moreover if also $w-2n\neq
0$, then
$$\operatorname{dim}_{\Q_p}H^1(K_{\p},M_{\theta\Q_p}(n))=[K_{\p}:\Q_p]\operatorname{dim}_{\Q_p}M_{\theta\Q_p}(n).$$
\end{cor}

\end{section}
\begin{section}{The conjecture in the critical situation}
We follow the notations of \S 2, but now $p$ does not necessarily
split in $K$. The critical situation corresponds to realizations
of the motive
$$M_{\theta}(n)$$
where $\theta$ has infinity type $(a,b)$ and $n$ satisfies:
$$\operatorname{min}(a,b)<n\leq \operatorname{max}(a,b).$$
Let us impose that the weight is $\leq -3$; this means $a+b-2n\leq
-3$.

Tsuji in \cite[\S 9,\S 10]{Tsu} proves the Jannsen conjecture in
this case, $p$ inert and $\theta$ of infinite type $(k-j,0)$ with
$\vartheta_{\theta}=(\lambda,\ldots,\lambda)$, $n=k$ and
$0\leq-j<k$. One obtains also the case $p$ inert and $\theta$ is
of infinite type $(0,j-k)$, $n=j$ and $0\leq -k<j$ by complex
conjugation. The $p$-adic realization corresponds to:
$$M_{(k-j,0)\Q_p}(k)\cong V_p(E)^{\otimes k}\otimes \overline{V_p(E)}^{\otimes j};$$
where $\otimes$ is $\otimes_{\Y_p\otimes_{\Y_K}K}$,
$\Y_p=\Z_p\otimes_{\Z}\Y_K$ and the bar means complex
conjugation. \\
Observe then that any $M_{\theta}(n)$ in the critical situation
has a $p$-adic realization (we write the situation
$b=\operatorname{min}(a,b)$) isomorphic to
$$V_p(E)^{\otimes n-b}\otimes \overline{V_p(E)}^{\otimes n-a}.$$

Thus, if we write $k=n-b$ and $j=n-a$ we need to study Jannsen's
conjecture only for the $p$-adic realizations $M_{(k-j,0)}(k)$
(the situation $a=\operatorname{min}(a,b)$ corresponds by complex
conjugation to study the motives $M_{(0,j-k)}(j)$).

We have the following isomorphism,
$$(M_{\theta\Q_p}(n))^*(1)\cong M_{\overline{\theta}\Q_p}(-n+w+1).$$
 After all the above considerations we have in the critical
 situation:

\begin{teo}[Tsuji, theorem 9.1 and prop.10.2\cite{Tsu}]\label{glotsu} Let $p$ be a prime $\geq 5$
which is inert in $K$ and such that $E$ has good reduction at the
prime above $p$ in $K$. Consider $M_{\theta}(n)$ with
$\operatorname{min}(a,b)<n\leq \operatorname{max}(a,b)$ ($a\neq
b$) satisfying $a+b-2n\leq -3$, then
\begin{enumerate}
\item $H^2(G_S,M_{\theta\Q_p}(n))=0$, and \item
$H^2(G_S,M_{\overline{\theta}\Q_p}(1-n+w))=0$, where
$G_S=\operatorname{Gal}(K_S/K)$ with $S$ a finite set of primes of
$K$ containing the primes dividing $p\f$.
\end{enumerate}
\end{teo}
\begin{teo}[Tsuji, lemma 10.1\cite{Tsu}] Let $\operatorname{min}(a,b)<n\leq
\operatorname{max}(a,b)$ ($a\neq b$) and $a+b-2n\leq -3$. Suppose
$p$ is inert in $K$ and let $\p$ be the prime of $K$ above $p$.
\begin{enumerate}
\item Let $\mathfrak{q}$ be any non-zero prime ideal of $\Y_K$
different to $\p$. Then,\\
$H^i(K_{\mathfrak{q}},M_{\theta\Q_p}(n))=0$ and
$H^i(K_{\mathfrak{q}},M_{\overline{\theta}\Q_p}(1-n+w))=0$ for all
$i$. \item $H^i(K_{\p},M_{\theta\Q_p}(n))=0$ and
$H^i(K_{\p},M_{\overline{\theta}\Q_p}(1-n+w))=0$ for $i=0$ and
$2$.
\end{enumerate}
\end{teo}
\begin{obs}\label{trestres} The main points needed to obtain
theorem \ref{glotsu} are \cite[Theorem 6.1]{Tsu} and the main
Iwasawa conjecture proved by Rubin. Both results are also known in
the split situation by B.Han \cite[\S 5.1]{Han} and Rubin
\cite{Ru3} respectively. Suppose now that $p$ splits in $K$,
$p=\p{\p}^*$. In order to obtain \ref{glotsu} for the split case
one could try to write in detail the second part of Tsuji's paper
\cite[II]{Tsu} replacing $\p$ by $p=\p{\p}^*$. Let us indicate
some steps:  rewrite \cite[Theorem 4\S 5.1]{Han} with the unit
$e_{1,S,p}$ in the notation of \cite[\S5]{Tsu} to obtain
\cite[Theorem 6.1]{Tsu} (one will need a result similar to
\cite[Lemma 3.3]{Ba2}), replace $\Y$ by $\Y_p=\Y_K\otimes\Z_p$ in
\S8 and \S9 of \cite{Tsu} and check that the arguments follow
(remember that in our case $L=K$ and
$A_{\mathfrak{p}}=A_{\mathfrak{\beta}}=\Y_p$).
\end{obs}
\end{section}
\begin{section}{
The conjecture in the non-critical situation}

Let us recall that we have fixed an imaginary quadratic field $K$
with $cl(K)=1$ and $E$ has good reduction at the places above $p$
in $K$.

>From the specialization of the elliptic polylogarithm \cite{K},
one can prove the equality between the image by the Soul\'e
regulator map $r_p$ of a module
$\mathcal{R}_{\theta}=\alpha_{\theta}\Y_K$ in $K$-theory for the
motive $e_{\theta}(\otimes^wh^1(E))(w+l+1)$ (see \cite[3.3]{Ba3}
for the precise definition of $\alpha_{\theta}$) and the image of
the elliptic units module by a Soul\'e map $e_p$ (see
\cite[4.6]{Ba3} for the definition of $e_p$) with
$w-2(w+l+1)\leq-3$ with $-l\leq \operatorname{min}(a,b)$ (we refer
to \cite[cor.5.9]{Ba3} for this equality). This result is
obtained, under some restrictions, from the main conjecture of
Iwasawa theory for imaginary quadratic fields.

 Let us fix in this section the
following assumptions: let $p\geq 5$ be a prime with
$p>N_{K/\Q}\f$ and $\theta$ with infinite type $(a,b)$ such that
$a\not\equiv b(mod|\Y_K^*|)$. Suppose also that $\Y_K^*\rightarrow
(\Y_K/\f_{\theta})^*$ is injective and the
$\Delta=\operatorname{Gal}(K(E[p])/K)$-representation of
$\operatorname{Hom}_{\Y_p}(M_{\theta\Z_p}(w+l),\Y_p)$ satisfies
the Rubin's theorem on the Iwasawa main conjecture \cite{Ru3},
where $\Y_p$ denotes $\Y_K\otimes_{\Z}\Z_p$.
\begin{prop} Suppose that
$H^2(\operatorname{Gal}(K_S/K),M_{\theta\Q_p}(w+l+1))=0$ with $S$
finite set of primes of $K$ which contains the ones dividing
$p\f$. Then $r_p(\mathcal{R}_{\theta}\otimes_{\Z}\Q_p)$ has
$\Q_p$-rank equal to 2, in particular $e_p$ is injective.
\end{prop}
\begin{dem} The same arguments of the proof of
\cite[Prop.5.2.5]{K} can be used with the diagram \cite[Lemma
4.11]{Ba3} to obtain the result.
\end{dem}
As observed by Tsuji \cite[(11.12)]{Tsu}, Jannsen's conjecture can
be obtained from the non-vanishing of the map $e_p$ into local
cohomology of the elliptic units module if \cite[Question
11.15]{Tsu} has a positive answer. In \cite{Ba3} we prove that
\cite[Question 11.15]{Tsu} has a positive answer in the split
case, and in this case one can show the same result of Tsuji by
using the arguments of remark \ref{trestres}. But since we do not
want to write all the technical details, we make the following

\vspace{0.4cm}

{\bf assumption:} {\bf (I)} the result \cite[(11.12)]{Tsu} is true
also when $p$ splits in $K$.

\vspace{0.4cm}

In his thesis \cite{G}, Geisser studies this Soul\'e map $e_p$
locally, i.e., with values in the local cohomology. In this study
appear the values at non-critical values of the $p$-adic
$L$-functions associated to $\psi_{\Omega_i}$ for $i=1,2$ (up to
twist by cyclotomic character). The $\psi_{\Omega_i}$ comes from
the Hecke character $\psi_{\theta}$ associated to $\theta$ as
follows:
$\psi_{\theta}\otimes\Z_p=\psi_{\Omega_1}\oplus\psi_{\Omega_2}$
(see for an extended explanation \cite{G} or \cite[\S3]{Ba2}). Let
us denote by
$\iota:H^1(\operatorname{Gal}(\overline{K}/K),M_{\theta\Z_p}(w+l+1))\rightarrow
H^1(\operatorname{Gal}(\overline{K_{\p}}/K_{\p}),M_{\theta\Z_p}(w+l+1))$
the restriction map.

\begin{prop}[3.4 in \cite{Ba2}]\label{archiv} Let $p$ be a prime that splits in
$K$.

Suppose furthermore $w=a+b\geq1$, $a+l>0$, $b+l>0$ and
$p>3w+2l+w+1$. Then the length of the coimage of $\iota\circ
r_p(\Rr_{\theta})$ in $H^1(K_{\p},M_{\theta\Z_p}(w+l+1))$ is equal
to the $p$-adic valuation of
$$G(\psi_{\Omega_1}\kappa^l,u_1^{-a_{\theta}-1}-1,u_2^{-b_{\theta}-1}-1)
G(\psi_{\Omega_2}\kappa^l,u_1^{-b_{\theta}-1}-1,u_2^{-a_{\theta}-1}-1),$$
where $\kappa$ denotes the cyclotomic character of
$\mathcal{G}=\operatorname{Gal}(K(E[p^{\infty}])/K)$ and $G$
denote the $p$-adic $L$-functions.
\end{prop}
\begin{cor}\label{tresdos} Let us assume (I) and the hypothesis of proposition \ref{archiv}.
If the non-critical values obtained in the $p$-adic $L$-functions
$G(\psi_{\Omega_1}\kappa^l,u_1^{-a_{\theta}-1}-1,u_2^{-b_{\theta}-1}-1)$
and
$G(\psi_{\Omega_2}\kappa^l,u_1^{-b_{\theta}-1}-1,u_2^{-a_{\theta}-1}-1)$
do not vanish, then Jannsen's conjecture is true for
$M_{\overline{\theta}}(-l)$, i.e.
$$H^2(G_S,M_{\overline{\theta}\Q_p}(-l))=0,$$
and the conclusions with $n=w+l+1$ of corollary \ref{dosdeu} are
also satisfied.
\end{cor}
\begin{dem}
By the local Jannsen conjecture we obtain from \ref{dosonze} that
$$rank_{\Z_p}H^1(K_{\p},M_{\theta\Z_p}(w+l+1))=2.$$ By the
non-vanishing hypotheses and \ref{archiv}, we see that the module
$\iota(r_p(\mathcal{R}_{\theta}))$ has finite index in
$H^1(K_{\p},M_{\theta\Z_p}(w+l+1))$; therefore
$\operatorname{dim}_{\Q_p}\iota(r_p(\mathcal{R}_{\theta}\otimes\Q_p))=\operatorname{dim}_{\Q_p}(\mathcal{R}_{\theta}\otimes\Q_p)$.
Hence
$\operatorname{dim}_{\Q_p}(r_p(\mathcal{R}_{\theta}\otimes\Q_p))=2$
and the Soul\'e map $r_p$ does not vanish in
$\mathcal{R}_{\theta}$. This is equivalent by \cite[(11.12)]{Tsu}
to Jannsen's conjecture for
 $H^2(G_S,M_{\overline{\theta}\Q_p}(-l))$ (using the Tate-Poitou long exact sequence).
\end{dem}
\begin{obs2}Moreover by \cite[p.171,(11.12)]{Tsu} (assuming (I)),
the Jannsen conjecture for $M_{\overline{\theta}\Q_p}(-l)$ implies
the non-vanishing of the product of the non-critical values of the
$p$-adic $L$-functions which appear in \ref{tresdos}.
\end{obs2}

\end{section}
\begin{section}{The weak Jannsen conjecture}
We want to consider now the following weak form of Jannsen's
conjecture, which claims for our concrete realizations that: the
Galois cohomology groups $H^2(G_S,M_{\theta\Q_p}(n))$ vanish for
almost all Tate twist $n$. This weak Jannsen conjecture is
equivalent to the weak Leopoldt conjecture \cite{PR1}.
\begin{prop}\label{wlcp} Let us suppose that the cyclotomic $\mu$-invariant of every
abelian extension of $K$ is zero. Let $S$ be a set of primes of
$K$ containing the primes of $K$ dividing $p\mathfrak{f}$.
 Then for almost all $n$ we have
 $$H^2(G_S,M_{\theta\Q_p}(n))=0.$$
\end{prop}
\begin{dem} Since we have a Hecke character the Galois group $\operatorname{Gal}(K'/K)$ fixing the
$p$-torsion of $M_{\theta\Q_p}(n)/M_{\theta\Z_p}(n)$ is an abelian
extension of $K$. By \cite[B.2 Corollary (i)]{PR1} (here we use
the hypothesis that $\mu$ vanishes for the cyclotomic extension of
the abelian field $K'$ over $K$) we obtain that
$$e_{\tilde{\chi}}H^2(\operatorname{Gal}(K_S/K(\mu_{p^{\infty}})),M_{\theta\Z_p}(n))$$ is a
$\Z_p[[\operatorname{Gal}(K(\mu_{p^{\infty}})/K)]]^{\tilde{\chi}}$-torsion
module, where $e_{\tilde{\chi}}$ are the idempotents associated to
characters $\tilde{\chi}$ of
$\operatorname{Gal}(K(\mu_{p^{\infty}})/C_{\infty})$, with
$C_{\infty}$ the cyclotomic extension of $K$. This fact is
equivalent to
$$H^2(\operatorname{Gal}(K_S/K(\mu_{p^{\infty}})),M_{\theta\Q_p}^*/M_{\theta\Z_p}^*)=0,$$
where $^*$ means $\operatorname{Hom}(,\Q_p)$ or
$\operatorname{Hom}(,\Z_p)$ respectively \cite[prop. 1.3.2]{PR1}.
Now use \cite[Lemma 8]{Ja} to obtain the result.
\end{dem}
\begin{obs2} The work of Kato
for CM modular forms \cite[\S 15]{Ka} gives a proof of the weak
Jannsen conjecture \ref{wlcp} without any assumption on the $\mu$
invariant. Moreover in \cite{Ba4} we obtain a different proof of
proposition \ref{wlcp} without the $\mu$ vanishing assumption
using Iwasawa modules in two variables, imposing that $p+1\nmid
|b-a|=a-b$ if $p$ is inert in $K$, or $p-1\nmid |b-a|=a-b$ if $p$
splits.
\end{obs2}
\end{section}
\begin{section}*{Acknowledgments}
I am very grateful to the referee for pointing out some
inaccuracies, for providing some ideas to prove the results in a
more clear way and for his/her useful comments and suggestions. I
am also very pleased to thank Xavier Xarles for his suggestions
and comments.
\end{section}

\vspace{0.3cm} Francesc Bars Cortina\\
\tiny{Departament de Matem\`atiques}\\
\tiny{Universitat Aut\`onoma de Barcelona}\\
\tiny{  08193 Bellaterra (Barcelona)}\\ \tiny{Catalunya, Spain}\\
\tiny{ e-mail: francesc@mat.uab.es}

\end{document}